\documentclass[12pt]{article}
\usepackage[russian,english]{babel}
\usepackage{amsfonts,amsmath,amssymb,amscd,amsthm,url}
\long\def\comment#1\endcomment{}

\def\sgn{\mathop{\fam0 sgn}}

    \theoremstyle{theorem}
         \newtheorem{Theorem}{Theorem}
         \newtheorem{Lemma}[Theorem]{Lemma}

\newtheoremstyle{mydefinition}
  {3pt}
  {3pt}
  {\normalfont}
  {\parindent}
  {\bfseries}
  {.}
  { }
  {}

\theoremstyle{mydefinition}

\begin{document}


\selectlanguage{english}

\centerline{\uppercase{\bf How Fermat found extrema}
\footnote{I am grateful to A. Doledenok, M. Skopenkov and A. Sgibnev for useful discussions.}
}

\bigskip
\centerline{\bf A. Skopenkov \footnote{\url{www.mccme.ru/~skopenko}.
Moscow Institute of Physics and Technology and Independent University of Moscow.
Supported in part by the D. Zimin Dynasty Foundation and Simons-IUM fellowship.
}}

\bigskip
{\bf Abstract.}
We present a short elementary proof of the well-known criterion for a cubic polynomial to have three real roots.
The proof is based on Fermat's approach to calculus for polynomials. 
This approach illustrates the idea of a derivative rigorously but without technical $\varepsilon$-$\delta$ language.
The note is accessible to high-school students.



\bigskip
In elementary school one learns the following.

{\it Let $a$ and $b$ be real numbers.
Then the following conditions are equivalent:

$\bullet$ there are real numbers $x,y$ such that $a=x+y$ and $b=xy$;

$\bullet$ the equation $t^2-at+b=0$ has a real root;

$\bullet$ $4b-a^2\le0.$}

In this note we illustrate the idea of a derivative by proving a well-known generalization to {\it three} numbers
(which might be proved by Fermat or earlier).

\begin{Theorem}\label{t:fermat}
Let $a,b,c$ be real numbers.
There are real numbers $x,y,z$ such that
$$a=x+y+z,\quad b=xy+yz+zx\quad\text{and}\quad c=xyz$$
if and only if
$$4\left(b-\dfrac{a^2}3\right)^3+27\left(c-\dfrac{ab}3+\dfrac{2a^3}{27}\right)^2\le0.$$
\end{Theorem}

This result is very useful; for recent applications to elementary inequalities see \cite{DFMS} and references therein.
Not only this result, but simple exposition of its proof is hopefully interesting.
We do not use the notion of a derivative.
However, we illustrate the idea of a derivative rigorously but without technical $\varepsilon$-$\delta$ language, cf. \cite[\S8.1, \S8.2]{ZSS}.
So this exposition might be useful for introductory courses on analysis.
Unfortunately, this exposition is not so well-known, cf. \cite{Po84}.\footnote{E.g. in June 2016 it was suggested to add this proof to \cite{DFMS} in order to make Theorem \ref{t:fermat}
and its applications more accessible.
However, the authors found this proof `too complicated'. 
The first version of this note (\url{https://arxiv.org/pdf/1610.05968v1.pdf}) is rejected from posting on \url{http://www.turgor.ru/lktg/2016} by the authors of \cite{DFMS} and S. Dorichenko.}
For development of the `graphs of functions' idea see \cite{FT, Go10, Ta88} (in particular, one can apply this idea to the $pqr$-lemmas of \cite{DFMS}).

There is an alternative proof of Theorem \ref{t:fermat} using complex numbers and calculation of the {\it discriminant} of cubic polynomial in terms of coefficients, see e.g. \cite[\S2 and solutions of problems 6-22]{DFMS}.
Although that proof is longer, it illustrates another interesting important ideas.
In \cite{Ta88} a geometric interpretation (but not a proof) of Theorem \ref{t:fermat} is presented.

\medskip
{\it A reformulation of Theorem \ref{t:fermat} in terms of cubic polynomials.}
By Vieta Theorem (i.e. by direct calculation) the condition on $x,y,z$ is equivalent to the following equality of polynomials:
$$t^3-at^2+bt-c=(t-x)(t-y)(t-z).$$

{\it A reduction of Theorem \ref{t:fermat} to the case $a=0$.}
Take $u:=t-\frac a3$. Then $t=u+\frac a3$, so
$$t^3-at^2+bt-c=u^3+\left(b-\dfrac{a^2}3\right)u-\left(c-\dfrac{ab}3+\dfrac{2a^3}{27}\right).$$
Thus it suffices to prove Theorem \ref{t:fermat} for $a=0$.

\begin{proof}[Proof of Theorem \ref{t:fermat} for $a=0$, $b\ge0$]
Assume that the required $x,y,z$ exist.
Since $b\ge0$, the function $t^3+bt-c$ is strictly increasing as a sum of increasing functions, one of them strictly increasing.
Then the equation $t^3+bt-c=0$ has at most one real root.
This and $a=0$ imply that $x=y=z=0$, hence $4b^3+27c^2=0$.

If $4b^3+27c^2\le0$, then $b=c=0$, so we can take $x=y=z=0$.
\end{proof}

{\it Heuristic considerations for investigation of the function $f(t):=t^3-3t$} (formally they are not used in the proof).
The function $f$ is called {\it strictly increasing} on an interval if $f(t_1)>f(t_2)$ for each different $t_1>t_2$ from the interval.
This is equivalent to $\varphi(t_1,t_2)>0$ for each different $t_1,t_2$ from the interval, where
$$(*)\qquad\varphi(t_1,t_2):=\dfrac{f(t_1)-f(t_2)}{t_1-t_2}=t_1^2+t_1t_2+t_2^2-3.$$
If these conditions hold for each two `close' $t_1,t_2$, then by transitivity they hold for each two $t_1,t_2$.
So we guess that the `boundary points' of the intervals on which $f$ is monotonous are the roots of the equation $t^2+tt+t^2-3=0$.
These roots are $\pm1$.
(This is analogous to \cite{Be88}; considering a simple example before general method makes the method more accessible.)

\begin{Lemma}\label{l:incr}
The function $f(t):=t^3-3t$ is strictly increasing on $(-\infty,-1]$, is strictly decreasing on $[-1,1]$ and is strictly increasing on $[1,+\infty)$.
\end{Lemma}

\begin{proof} Denote $\varphi(t_1,t_2)$ by the formula (*).
Then $\varphi(t_1,t_2)>0$ for each different $t_1,t_2\ge1$.
Hence $f(t)$ strictly increases on $[1,+\infty)$.
Analogously the other two statements are proved.
\end{proof}

\begin{proof}[Proof of Theorem \ref{t:fermat} for $a=0$, $b<0$]
Taking
$u:=t\sqrt{-b/3}$ we may assume that $b=-3$.

First assume that the required $x,y,z$ exist.
Since $-3=b<0=a$, the case $x=y=z$ is impossible, so the equation $f(t)=t^3-3t-c=0$ has at least two real roots.

Denote $f(t):=t^3-3t-c=t(t^2-3)-c$.
Then $t_+:=1+\max\{3,|c|\}>1$ and $f(t_+)>0$.
Analogously, there is $t_-<-1$ such that $f(t_-)<0$.
So by Lemma \ref{l:incr} and the Intermediate Value Theorem, the equation $f(t)=0$ has at least two real roots if and only if $f(-1)$ and $f(1)$ `have different signs', i.e. $f(-1)f(1)\le0$.
We have $-f(t)=c-t(t^2-3)$, so
$$f(-1)f(1)=(c+2\cdot1)(c-2\cdot1)=c^2-4=(4b^3+27c^2)/27.$$

Now assume that $4b^3+27c^2=27(c^2-4)\le0$.
If $c^2<4$, the above paragraph shows that the equation $f(t)=0$ has three real roots, and we take them as $x,y,z$.
If $c^2=4$, we take $x=y=-\sgn c$ and $z=2\sgn c$ (it  should be clear how to guess these formulas).
\end{proof}

{\bf Exercise.} Find maximal intervals on which the following function is strictly increasing (decreasing):

(a) $f(t)=t^4-4t$. \quad (b) $f(t)=t^4-12t^3+22t^2-24t+10$.


\end{document}


\title{Как Ферма искал экстремумы}


\author{А. Скопенков}

\thanks{Благодарю М. Горелова, А. Доледенка, М. Скопенкова, А. Сгибнева и анонимного рецензента за полезные замечания.
\newline
\url{www.mccme.ru/~skopenko}.
Московский физико-технический институт и Независимый Московский Университет.
Частично поддержан грантами фонда Д. Зимина Династия и фонда Саймонса.}

\date{}

\maketitle

\abstract
Приведено элементарное доказательство критерия существования трех различных вещественных корней у многочлена третьей степени.
В основу доказательства лег подход Ферма к анализу для многочленов.
Этот подход иллюстрирует понятие производной на математическом уровне строгости, но без использования языка $\varepsilon$-$\delta$.
Статья доступна для понимания старшеклассников.
\endabstract

\section{О чем эта статья}

В этой статье мы на конкретном примере предложим решение одной из важнейших проблем преподавания математичексого анализа --- как показать его важнейшие идеи на математическом уровне строгости, но без 
использования технического языка $\varepsilon$-$\delta$.\footnote{Для нахождения количества корней (но не для нахождения экстремумов) мы используем теорему \ref{t:inter} о промежуточном значении.
Эта теорема интуитивно очевидна; ее доказательство необходимо при {\it любом} подходе к нахождению количества корней кубического трехчлена.}
(Естественно, перед этим полезно знакомство с идеями на физическом уровне строгости.) 
Иными словами, мы выводим на первый план <<олимпиадную>> составляющую изучения математического анализа.
Поэтому приводимое изложение может быть полезно при изучении его основ.
К сожалению, оно малоизвестно.\footnote{Ср. \cite{Po84}. Другой пример: в июне 2016 было предложено добавить это доказательство в \cite{DFMS} с целью сделать теорему \ref{t:fermat} и ее применения более доступными.
Однако, авторы сочли доказательство <<слишком сложным>>. 
Первая версия этой статьи (\url{https://arxiv.org/pdf/1610.05968v1.pdf}) не принята к выкладыванию
на \url{http://www.turgor.ru/lktg/2016} авторами статьи \cite{DFMS} и С. Дориченко.}

Наш конкретный пример --- нахождение минимумов и максимумов многочлена, которое позволяет находить количество его корней (не находя самих корней).
Точнее, мы сведем задачу о нахождении экстремумов к задаче о поиске корней многочлена меньшей степени.

Приводимый подход принадлежит Пьеру Ферма (1601-1665);  теоремы \ref{t:fermat0s}, \ref{t:mono30}, \ref{t:mono3}, \ref{t:fermat0} и \ref{t:fermat} получены им или еще ранее.
Об увлекательной истории этих открытий написано, например, в \cite{Gi82, Pr18, St84, Yu70}.

О дальнейшем развитии идеи Ферма об <<анализе для многочленов>>  см. задачи \ref{p:gen}, \ref{numro4} и \cite{Lu14, Lu}, \cite[\S8.2, \S8.3]{ZSS} (полезная задача --- аналогично данной статье извлечь строгие доказательства из эвристических рассуждений в \cite{Lu}).
Развитие идеи <<графиков функций>> описано в \cite{FT, Go09, Go10, Ta88}
(например, эта идея может быть применена к  $pqr$-леммам в \cite{DFMS}).


Чтобы заинтересовать читателя приводимым методом (\S\ref{s:analysis}), перед его изложением мы приведем некоторые результаты, которые могут быть им получены (теоремы \ref{t:fermat0s}, \ref{t:mono30}, \ref{t:mono3}, \ref{t:fermat0} и \ref{t:fermat}).
Такая структура текста покажет, что можно легко применять результат, не вникая в метод его придумывания и доказательства.
Думаю, тогда читателю будет интересно вникать в этот метод.
Именно он (а не приводимые результаты) составляет основное содержание данной заметки.

Читателю, не знакомому с математическим анализом, удобнее иметь прямые элементарные формулировку и доказательство результата, а не вывод результата из менее элементарной формулировки, предваренной немотивированной теорией.
Надеюсь, это и мотивирует его к изучению теории (коль скоро его заинтересуют обобщения), и поможет ему изучить теорию.

\section{Формулировки}

В школе изучают следующие факты.

\begin{Theorem}\label{t:fermat2}
Пусть $a$ и $b$ --- вещественные числа.
Тогда следующие условия эквивалентны:

(1) существуют вещественные числа $x,y$ такие, что $a=x+y$ и $b=xy$;

(2) уравнение $t^2-at+b=0$ имеет вещественный корень;

(3) $4b-a^2\le0.$
\end{Theorem}

Вот близкий результат: квадратное уравнение $t^2-at+b=0$ имеет два решения при $D:=a^2-4b>0$, имеет одно решение при $D=0$ и~не имеет решений при $D<0$.

Заметим, что модуль дискриминанта $D$ --- это расстояние от вершины параболы до оси абсцисс.

Функция $f$ называется {\bf строго возрастающей} на интервале, если $f(t_1)>f(t_2)$ для любых $t_1>t_2$ из интервала.
Аналогично определяется строгое убывание.

\begin{Theorem}\label{t:mono2}
Пусть $a$ и $b$ --- вещественные числа.
Тогда функция $t^2-at+b$ строго убывает на $(-\infty,a/2]$ и строго возрастает на $[a/2,+\infty)$.
\end{Theorem}

Эти теоремы (как и формула для корней квадратного уравнения) доказываются при помощи равенства
$$t^2-at+b=\left(t-\frac a2\right)^2+\left(b-\frac{a^2}4\right).$$
Мы проиллюстрируем понятие производной на примере доказательства хорошо известного обобщения
приведенных теорем на {\it три} числа.
В этом параграфе мы приводим <<ведущие>> частные случаи.
Формулировки результатов для общего случая (теорем \ref{t:fermat} и \ref{t:mono3}) немного более громоздкие.
А их доказательства состоят в несложной редукции к <<ведущим>> частным случаям.
Они приводится в \S\ref{s:algebre} и \S\ref{s:analysis}.


\begin{Theorem}\label{t:fermat0s}
Пусть $b$ и $c$ --- вещественные числа.
Тогда следующие условия эквивалентны:

(1) существуют попарно различные вещественные числа $x,y,z$, такие что
$$0=x+y+z,\quad b=xy+yz+zx\quad\text{и}\quad c=xyz.$$

(2) уравнение $t^3+bt-c=0$ имеет три попарно различных вещественных корня.

(3) $4b^3+27c^2<0$.
\end{Theorem}

Этот результат (точнее, обобщающая его нижеприведенная теорема \ref{t:fermat0}) очень полезен.
О применениях к элементарным неравенствам см., например, \cite[задачи 13 и 32]{SB78}, \cite{Go09, Go}; \cite[задача 2]{DFMS} --- частный случай равносильной версии условия (3) теоремы \ref{t:fermat0s} (точнее,  аналогичной равносильной версии условия (3) нижеприведенной теоремы \ref{t:fermat0}).

Заметим, что

$\bullet$ условие (3) заведомо не выполнено при $b\ge0$,

$\bullet$ условие (3) равносильно условию $b<-3\sqrt[3]{c^2/4}$ или условию $|c|<2\sqrt{-b^3/27}$,

\begin{figure}[h]\centering
\centerline{\includegraphics{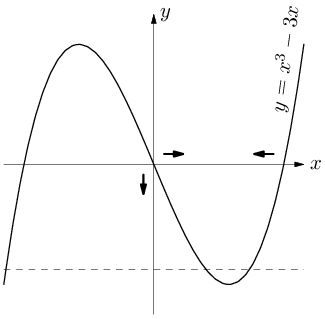}}
\caption{Графики функций $f(t)=t^3+t$ и $f(t)=t^3-3t$ (нужно отредактировать рисунок и добавить график функции $f(t)=t^3+t$)} \label{f:b}
\end{figure}


$\bullet$ выражение $4b^3+27c^2$ пропорционально произведению ординат точек экстремума функции $f(t)=t^3+bt-c$ (см. следующие лемму и теорему).

\begin{Lemma}\label{l:discr} Для любых $b\le0$ и $c$ обозначим $f(t):=t^3+bt-c$.
Тогда $f(\sqrt{-b/3})f(-\sqrt{-b/3})=\frac{4b^3}{27}+c^2$.
\end{Lemma}

\begin{Theorem}\label{t:mono30}
Пусть $b$ и $c$ --- вещественные числа.
Тогда функция $t^3+bt-c$

$\bullet$ при $b\ge0$ строго возрастает на $(-\infty,+\infty)$;

$\bullet$ при $b<0$ строго возрастает на $(-\infty,-\sqrt{-b/3}]$, строго убывает на $[-\sqrt{-b/3},\sqrt{-b/3}]$ и строго возрастает на $[\sqrt{-b/3},+\infty)$.
\end{Theorem}

Заметим, что метод Ферма иллюстрируется именно на доказательстве теоремы \ref{t:mono30}.
Она применяется для доказательства наиболее нетривиальной части теоремы \ref{t:fermat0s}.


\section{Анализ: доказательство теоремы \ref{t:mono30}}\label{s:analysis}


{\it Эвристические соображения для вывода теоремы \ref{t:mono30}} (формально, они не используются в доказательстве).
Свободный член $c$ при исследовании на возрастание и убывание не играет роли.
На рис. \ref{f:b} приведены графики функции $f(t)=t^3+bt$ при разных $b$.
Ясно, что

$\bullet$ при $b\ge0$ функция $f(t)$ возрастает, и

$\bullet$ при $b<0$ функция $f(t)$ имеет локальный максимум и локальный минимум.

Покажем, как найти  локальный максимум и локальный минимум на примере случая $b=-3$ (общий случай можно свести к нему, что сделано ниже, или рассмотреть аналогично).
Т.е. найдем локальный максимум и локальный минимум функции $f(t)=t^3-3t$.
Условие строгого возрастания функции $f$ равносильно условию $\varphi(t_1,t_2)>0$ для любых различных $t_1,t_2$, где
$$(*)\qquad\varphi(t_1,t_2):= \dfrac{f(t_1)-f(t_2)}{t_1-t_2}= \dfrac{t_1^3-t_2^3-3(t_1-t_2)}{t_1-t_2}= t_1^2+t_1t_2+t_2^2-3.$$
Если эти условия выполнены для любых двух <<близких по значению>> $t_1,t_2$,
то в силу транзитивности неравенств они выполняются для любых двух $t_1,t_2$.
Так мы приходим к догадке, что граничные точки интервалов, на которых $f$ монотонна, являются корнями уравнения $t^2+tt+t^2-3=0$.
Эти корни равны $\pm1$.
(Это рассуждение аналогично \cite{Be88}; разбор простого примера перед рассказом общего метода делает этот метод более понятным.)

Смысл этого построения --- замена сложного условия $\varphi(t_1,t_2)>0$ возрастания функции на гораздо более простое условие $\varphi(t,t)>0$ (в нашем случае это просто $3t^2-3>0$).
Гениальная догадка первопроходцев математического анализа состояла в том, что
<<в большинстве интересных случаев>>  эти условия равносильны.
Эта догадка многократно использовалась Ферма, Ньютоном, Лейбницем, Эйлером и другими.
Четкие определения и строгие доказательства появились только в XIX веке.
Функция $\varphi(t,t)$ называется {\it производной} функции $f$.


\begin{proof}[Доказательство теоремы \ref{t:mono30}]
Функция $t^3$ строго возрастающая (это доказывается перемножением неравенств).
Аналогично при $b\ge0$ функция $bt$ возрастающая.
Поэтому при $b\ge0$ функция $t^3+bt-c$ строго возрастающая.

Пусть теперь $b<0$.
Можно считать, что $c=0$.
Будем считать, что $b=-3$ (см. начало доказательства леммы \ref{l:discr} в \ref{s:algebre}).
Введем $\varphi(t_1,t_2)$ по формуле (*) для $f(t):=t^3-3t$.
Тогда $\varphi(t_1,t_2)>0$ для любых различных $t_1,t_2\ge1$.
Следовательно, $f(t)$ строго возрастает на $[1,+\infty)$.
Два другие утверждения теоремы доказываются аналогично.
\end{proof}

\begin{Theorem}\label{t:mono3}
Пусть $a,b,c$ --- вещественные числа.
Тогда функция $t^3-at^2+bt-c$

$\bullet$ при $3b\ge a^2$ строго возрастает на $(-\infty,+\infty)$;

$\bullet$ при $3b<a^2$ строго возрастает на $\left(-\infty,\dfrac{a-\delta}3\right]$, строго убывает на
$\left[\dfrac{a-\delta}3,\dfrac{a+\delta}3\right]$ и строго возрастает на $\left[\dfrac{a+\delta}3,+\infty\right)$, где $\delta=\sqrt{a^2-3b}$.
\end{Theorem}

Докажите эту теорему (например, сведя ее аналогично \S\ref{s:algebre} к частному случаю $a=0$, т.е. к теореме \ref{t:mono30})!

\begin{pr}\label{p:gen} Найдите максимальные интервалы строгого возрастания (убывания) функции

(a) $f(t)=t^4-4t$; \quad (b) $f(t)=t^4-12t^3+22t^2-24t+10$.
\end{pr}

\section{Алгебра: доказательства теорем \ref{t:fermat0s}, \ref{t:fermat0} и \ref{t:fermat}}\label{s:algebre}

Теорема \ref{t:fermat0s} доказывается аналогично следующему результату.

\begin{Theorem}\label{t:fermat0}
Пусть $b$ и $c$ --- вещественные числа.
Тогда следующие условия эквивалентны:

(1) существуют вещественные числа $x,y,z$, такие что
$$0=x+y+z,\quad b=xy+yz+zx\quad\text{и}\quad c=xyz.$$

(2) уравнение $t^3+bt-c=0$ имеет три вещественных корня с учетом кратности.


(3) $4b^3+27c^2\le0$.
\end{Theorem}

Наличие трех вещественных корней с учетом кратности у уравнения $t^3-at^2+bt-c=0$ определяется как наличие    вещественных чисел $x,y,z$, для которых равны многочлены
$$t^3-at^2+bt-c=(t-x)(t-y)(t-z)$$
(это равенство означает, что равны коэффициенты при соответствующих степенях).
Раскрывая скобки, видим, что это равносильно условию (1). 
(Заодно мы получили {\it теорему Виета} для кубического уравнения.)
Мы доказали равносильность $(1)\Leftrightarrow(2)$ сразу для следующего общего случая.

\begin{Theorem}\label{t:fermat}
Пусть $a,b,c$ --- вещественные числа.
Тогда следующие условия эквивалентны:

(1) существуют вещественные числа $x,y,z$, такие что
$$a=x+y+z,\quad b=xy+yz+zx\quad\text{и}\quad c=xyz.$$

(2) уравнение $t^3-at^2+bt-c=0$ имеет три вещественных корня с учетом кратности.

(3) $4\left(b-\dfrac{a^2}3\right)^3+27\left(c-\dfrac{ab}3+\dfrac{2a^3}{27}\right)^2\le0.$
\end{Theorem}

{\it Сведение теоремы \ref{t:fermat} к частному случаю $a=0$, т.е. к теореме \ref{t:fermat0}.}
Положим $u:=t-\frac a3$. Тогда $t=u+\frac a3$, поэтому
$$t^3-at^2+bt-c=u^3+\left(b-\dfrac{a^2}3\right)u-\left(c-\dfrac{ab}3+\dfrac{2a^3}{27}\right).$$
Следовательно, теорема \ref{t:fermat} следует из частного случая $a=0$, т.е. из теоремы \ref{t:fermat0}.

\medskip
Чтобы продемонстрировать полезность нахождения промежутков возрастания и убывания функций, мы выведем равносильность $(2)\Leftrightarrow(3)$ в теореме \ref{t:fermat0} из теоремы \ref{t:mono30}.
После этого приводится более простое доказательство этой равносильности (для нетривиального случая $b<0$), не использующее  теоремы \ref{t:mono30} и сообщенное мне М. Гореловым.
Имеется еще одно доказательство теорем \ref{t:fermat0} и \ref{t:fermat} с использованием комплексных чисел и вычисления {\it дискриминанта} кубического многочлена в терминах коэффициентов, см. например \cite[\S2 и решения задач 6-22]{DFMS}.
Хотя это доказательство длиннее каждого из приведенных, оно иллюстрирует другие важные интересные идеи.
В \cite{Ta88} представлена геометрическая интерпретация (но не доказательство) теорем \ref{t:fermat0} и \ref{t:fermat}.

\begin{proof}[Доказательство равносильности $(2)\Leftrightarrow(3)$ в теореме \ref{t:fermat0} для $b\ge0$]
Сначала предположим, что требуемые корни $x,y,z$ существуют.
Так как $b\ge0$, то функция $t^3+bt-c$ строго возрастающая (это тривиальный случай теоремы \ref{t:mono30}).
Поэтому уравнение $t^3+bt-c=0$ имеет не больше одного действительного корня.
Значит $x=y=z$. Этот факт и условие $x+y+z=0$ влекут равенство $x=y=z=0$.
Следовательно $4b^3+27c^2=0$.

Теперь предположим, что $4b^3+27c^2\le0$. Тогда $b=c=0$, поэтому можем взять $x=y=z=0$.
\end{proof}

\begin{proof}[Доказательство леммы \ref{l:discr}]
Будем считать, что $b=-3$ (общий случай сводится к этому частному заменой $u:=t\sqrt{-b/3}$).
Имеем $-f(t)=c-t(t^2-3)$, поэтому
$$f(-1)f(1)=(c+2\cdot1)(c-2\cdot1)=c^2-4=(4b^3+27c^2)/27.$$
\end{proof}

\begin{proof}[Доказательство равносильности $(2)\Leftrightarrow(3)$ в теореме \ref{t:fermat0} для $b<0$]
Будем считать, что $b=-3$ (см. начало доказательства леммы \ref{l:discr}).
Обозначим $f(t):=t^3-3t-c$.
По лемме \ref{l:discr} неравенство (3) равносильно условию $f(-1)f(1)\le0$.

Предположим, что требуемые в условии (2) корни $x,y,z$ существуют.
Так как $xy+yz+zx=b=-3<0$, то случай $x=y=z$ невозможен.
Поэтому уравнение $f(t)=0$ имеет по крайней мере два различных вещественных корня.
Тогда по теореме \ref{t:mono30} для $b=-3$ числа $f(-1)$ и $f(1)$ <<имеют различные знаки>>, то есть $f(-1)f(1)\le0$.
(Действительно, если $f(-1)>0$ и $f(1)>0$, то по теореме \ref{t:mono30} уравнение $f(t)=0$ имеет не более одного вещественного корня.
Аналогично разбирается случай, когда $f(-1)<0$ и $f(1)<0$.)

Теперь предположим, что $f(-1)f(1)\le0$.

Пусть сначала $f(-1)f(1)=0$.
Тогда $c^2-4=0$.
Положим $x=y=-\sgn c$ и $z=2\sgn c$ (читатель наверняка догадался, как подобрать эти формулы).

Пусть теперь $f(-1)f(1)<0$.
Обозначим через $t_+$ наибольшее из чисел 2 и $1+|c|$.
Тогда
$$t_+\ge 1+|c|>1\quad\text{и}\quad f(t_+)=t_+(t_+^2-3)-c>(1+|c|)(2^2-3)-c>0.$$
Аналогично доказывается, что существует $t_-<-1$ такое, что $f(t_-)<0$.
Так как $t_-<-1$, $t_+>1$, $f(t_-)<0<f(t_+)$ и $f(-1)f(1)<0$, то по теореме \ref{t:mono30} и теореме \ref{t:inter} о промежуточном значении (см. ниже) уравнение $f(t)=0$ имеет три вещественных корня.
Мы возьмем их в качестве $x,y,z$.
\end{proof}

\begin{Theorem}[о~промежуточном значении]\label{t:inter}
Для многочлена $f$ и~чисел $a<b$ если $f(a)>0>f(b)$, то существует такое $c\in(a,b)$, что $f(c) =0$.
\end{Theorem}

Строгое доказательство этой интуитивно очевидной теоремы несложно, но требует аккуратного построения теории
вещественных чисел.
Мы его не приводим, см. \cite{Zo15}; см. также \cite{BG}.


\medskip
{\it Набросок другого доказательства равносильности $(2)\Leftrightarrow(3)$ в теореме \ref{t:fermat0} для $b<0$.}
Сначала повторяем первый абзац предыдущего доказательства.

Предположим, что требуемые корни $x,y,z$ существуют.
Будем считать, что все они отличны от $\pm1$ (этот случай можно рассмотреть отдельно).
Тогда из равенства $-3=xy+xz+yz=-x^2-xy-y^2$ следует, что на каждом из интервалов $(-\infty,–1),(-1,1),(1,+\infty)$ находится не более одного корня.
Поэтому на каждом из интервалов есть ровно один корень.
Тогда из равенства $f(t)=(t-x)(t-y)(t-z)$ следует, что $f(-1)>0$ и $f(1)<0$.
Поэтому $f(-1)f(1)<0$.

Теперь предположим, что $f(-1)f(1)\le0$.
Тогда по теореме о промежуточном значении уравнение $f(t)=0$ имеет корень $x\in[-1,1]$.
Значит,
$$f(t)=f(t)-f(x)=(t-x)(t^2 +xt+(x^2-3)).$$
Будем считать, что $x\ne\pm1$ (этот случай можно рассмотреть отдельно).
Так как $x\in(-1,1)$, то дискриминант квадратного трехчлена $t^2 +xt+(x^2-3)$ (от $t$) положителен.
Значит, по теореме \ref{t:fermat2} уравнение $f(t)=0$ имеет еще два корня $y,z$.
\qed

\begin{pr}\label{numro4}
Как по числам $n,p,q$ определить количество решений уравнения

(a) $t^4-4t+q=0$; \quad (b) $t^4+pt+q=0$; \quad (c) $t^n+pt+q=0$?
\end{pr}
